\documentclass[12pt]{iopart}

\usepackage{amsfonts,amssymb,amsbsy}
 
\usepackage{color}

\newcommand{\redc}{}

\expandafter\let\csname equation*\endcsname\relax
\expandafter\let\csname endequation*\endcsname\relax
\usepackage{amsmath} 

\newtheorem{proposition}{Proposition}
\newtheorem{remark}{Remark}

\def\debproof{ {\bf Proof.} }
\def\finproof{\hfill {\small $\Box$} \\}  

\newcommand{\eps}{\varepsilon}
\newcommand{\EE}{\mathbb{E}}

\begin{document}

\title{Global acoustic daylight imaging in a stratified Earth-like model}

\author{Maarten V. de Hoop$^1$,
Josselin Garnier$^2$, and
Knut S\o lna$^3$}

\address{$^1$ Department of Computational and Applied Mathematics and Department of Earth Science, Rice University, 6100 Main Street, Houston TX 77005}
\address{$^2$ Centre de Math\'ematiques Appliqu\'ees, Ecole Polytechnique,
91128 Palaiseau Cedex, France}
\address{$^3$ Department of Mathematics, University at California
  at Irvine, Irvine, CA 92697}

\ead{mailto{mdehoop@rice.edu}, \mailto{josselin.garnier@polytechnique.edu}, \mailto{ksolna@math.uci.edu}}

\begin{abstract}
We present an analysis of acoustic daylight imaging in an Earth-like
model assuming a random distribution of noise sources spatially supported in
an annulus located away from the surface. 
We assume a situation with scalar wave propagation and that the  measurements are
of  the wave field at the surface.
Then, we  obtain a relation between
the autocorrelation function of the measurements and the trace of the
scattered field generated by an impulsive source localized just below
the surface. From this relation it is, for example, clear that the
eigenfrequencies can be recovered from the autocorrelation. 
Moreover, the complete scattering operator can be extracted under 
the additional assumption that the annulus is close to the surface and 
has a thickness  smaller than the typical wavelength.
\end{abstract}


\submitto{\IP}

\maketitle

%\pagestyle{myheadings}
%\thispagestyle{plain}
%\markboth{M. de Hoop, J. Garnier, and K. S\o lna}{Global acoustic daylight imaging in a stratified
%       Earth-like model}

\section{Introduction} 
\label{intro}
The emergence of the Green's function from cross correlations of noisy
signals, in the context of seismic exploration, was first pointed out
by Claerbout \cite{claerbout68,claerbout85,claerbout99}. 
He modeled the  Earth's crust 
as a half-space and considered the situation
when waves are generated by  unknown sources in the
crust and recorded at the surface.
He showed  that the autocorrelation function of the
such signals  is the same as
the signal reflected by the Earth's crust  when an impulsive source is
used at the surface. 
The latter signal  is indeed the one recorded in reflection
seismology.  
He thereby established 
a formal connection between reflection seismology and wave field 
correlations.
The process of  autocorrelating signals  on
the surface, signals  that are  generated by noise sources in the interior,  
has come to be known as the  
daylight (imaging) configuration.  

The physical explanation of why daylight
imaging is equivalent to reflection seismology was simple and based on
flux conservation. It is possible to give a mathematical proof in the
case of a one-dimensional half space with radiation condition
\cite{garnier16,noisebook}.

Imaging of the structure of the Earth, from the core to the surface,
remains a major challenge. On the one hand, the standard seismic
methods using earthquake signals are limited by the uneven spatial
distribution of earthquakes, along the major tectonic plate
boundaries. \redc{On the other hand, seismic interferometry using ambient
noise records has been successful.
We here refer to seismic interferometry as using  cross correlations of signals 
to gain useful information about the medium.  
 Claerbout   in his book \cite{claerbout76}  asks the reader in an exercise 
to prove that the temporal autocorrelation of a transmission seismogram
 for a layered model with a source underground is equivalent to a reflection seismogram, a problem he treated in
 \cite{claerbout76} for a Goupillaud medium. 
Such a  
``daylight'' imaging approach has been subsequently validated experimentally and  
 mathematically for more general microstructures than layered, see for instance \cite{noisebook}. 
Seismic interferometry in general has  mostly been used to probe the uppermost
layers of the Earth through surface waves
\cite{bensen06,shapiro04}. However, recently, body-wave path
responses from 
 probing the deepest part of the Earth were obtained from noise records
\cite{boue13}. 
This motivates the question:  can the  daylight 
imaging approach as proposed by Claerbout  be extended 
to a global image procedure  for the Earth?
One  caveat is that in global seismology the model domain is a ball,
whence the radiation condition that played a crucial role in the
analysis of daylight imaging cannot be used anymore.  } If the noise
sources were distributed uniformly, then we would have equipartition
of energy amongst the normal modes, and we could invoke this argument
to establish the relation between the autocorrelation function of the
noise signals and the trace of the scattered field generated by
localized and impulsive sources just below the surface
\cite{weaver01}.  Unfortunately, the noise source distribution is not
uniform, and we cannot invoke the diffusion or ergodicity properties
of the Earth to claim that the observed field is equipartitioned as in
\cite{bardos08,garniersolna11} because the scattering is not strong enough. 
In  fact, we will show below  that in the propagation
regime of interest, the relation between the autocorrelation function
of the noise signals and the mentioned scattered field is not as
simple as in the case of  a one-dimensional half
space. 
In the case of a spherically symmetric Earth model, we show
that the point spectrum of the Earth  associated with the low angular orders 
can be extracted from  the correlation functions 
of the signals recorded at the surface and emitted by unknown noise sources localized 
away from the surface.
If  the noise sources satisfy some specific conditions in terms of
their spatial localization then  the correlations of noisy
signals can be used to extract the standard scattering operator. 
These conditions, that the source distribution is spatially localized 
on a thin annulus  with thickness comparable to or smaller than the typical
wavelength, could be satisfied in practice. We also include weak
angular variations of  the wavespeed in the analysis, which confirms 
that global acoustic  daylight imaging is possible
for an essentially spherically symmetric Earth.

The paper is organized as follows. 
\redc{In Section~2 we give a summary of the main results of the paper.
In Section~3 we introduce
the wavefield decomposition and coupling in a radial Earth and the
relevant asymptotic analysis. 
The scaling in our analysis corresponds to modes that essentially propagate 
in the radial direction.
In Section~4 we describe the scattering
operator in a stratified spherically symmetric Earth. 
In Section~5, we present global acoustic daylight imaging
and we establish the relationship between the scattering operator and
the  autocorrelation functions of the ambient noise signals.
We also consider the special case where
the spatial support of the noise sources is a thin annulus located below the surface 
with thickness smaller than the typical wavelength,
 providing a simplified result.  In Section~6, we  {show} 
 that our results are robust with respect to the properties of source and receiver distributions,
 with respect to small angular undulations  in the Earth's parameters 
(rather than a purely radial Earth), and with respect to measurement noise.
In particular,  it follows that even with small angular variations 
the observed point spectrum of the Earth can be reasonably accurately explained by
a radial model.}

\redc{\section{Summary of Modeling and Main Results} }
In this section we summarize the main modeling assumptions, quantities of interest,
and main result. 
\redc{We compare the symmetrized field (\ref{eq:symp1}) measured at the Earth's surface
and transmitted by an impulsive point source just below the surface with the empirical autocorrelation function (\ref{eq:CT1})
of the field measured at the Earth's surface and generated by a distributed noise source distribution.
We establish a correspondence  between them and we briefly discuss aspects of  robustness of the results to the modeling 
assumptions.  }

\subsection{Wave Decomposition in the Radial Earth and Main Result}
We consider here scalar wave propagation in the spherical Earth so that 
the harmonic wave field $\hat{p}= \hat{p}(\omega,r,\theta,\varphi)$ satisfies 
 \begin{eqnarray}\label{eq:main2}
\fl
\frac{\partial}{\partial r} r^2 \frac{\partial}{\partial r} \hat{p}  +\frac{1}{\sin \theta}
\frac{\partial}{\partial \theta} \sin \theta \frac{\partial}{\partial \theta} \hat{p} 
+\frac{1}{\sin^2 \theta} \frac{\partial^2 }{\partial \varphi^2} \hat{p} 
+ \frac{\omega^2 r^2}{  c^2(r)} \hat{p}  =   \hat{f} (\omega,r,\theta,\varphi) ,
\end{eqnarray}
for $(r,\theta,\varphi) \in (0,R_o) \times (0,\pi) \times(0,2\pi)$,
where $c(r)$ is the {\it radial} velocity model and $\hat{f}(\omega,r,\theta,\varphi)$
is the source term in the frequency domain. 
We assume here that 
$\hat{f}$ is supported away from the origin in $\omega$ and 
a boundary condition at the surface $r=R_o$ that takes 
into account boundary dissipation,  see Eq.~(\ref{eq:dissbc}).
We moreover assume a high-frequency situation, see Section \ref{sec:HF} for 
a discussion of the WKB context that this entails and which is central to our analysis.
 The wave field is then decomposed as
\begin{eqnarray}
\label{eq:fieldmode1a0}
\hat{p} (\omega,r,\theta,\varphi) =& \sum_{l,m} Y_{l,m}(\theta,\varphi)\hat{p}_{l,m} (\omega,r),  \end{eqnarray}
with the spherical harmonics $Y_{l,m}$ defined in Eq.  (\ref{eq:Ydef}) and 
\begin{eqnarray}
\label{eq:fieldmode1a00}
  \hat{p}_{l,m}(\omega,r) &=&  \int_0^\pi \int_0^{ 2\pi} \overline{Y_{l,m}(\theta,\varphi) }
 \hat{p}(\omega,r,\theta,\varphi) \sin \theta d\varphi d\theta  .
   %  \\      &=& F_{l,m} \delta(r-R_s) . 
\end{eqnarray}

We consider the two following source scenarios.

\subsubsection{Point Source.}
  
 This corresponds to a ``classic'' seismology configuration with a point  source
 just below the surface.
Let    $\hat{\cal R}_{l}(\omega,R_o)$ be the 
fundamental scattering function  associated  with Earth
with no surface reflection, note that this is independent of
$m$ as the velocity model is radial (see  Proposition \ref{prop1}).  
Then the symmetrized  surface Earth response function  is (see Proposition \ref{prop2})
\begin{eqnarray}
\label{eq:SO}
\hat{\cal S}_{l}(\omega,R_o) = 
 \left|
 \frac{1 + \hat{\cal R}_{l} (\omega,R_o)}
{1- \Gamma_{R_o} \hat{\cal R}_{l} (\omega,R_o)} \right|^2,
\end{eqnarray}
and this is the function that encapsulates information about the Earth's
interior. 
The 1  in the numerator corresponds to a direct transmission 
from the source and the $\hat{\cal R}_{l}$  term  is the reflection from the Earth's
interior 
while the denominator produces multiples from reflections at the Earth's surface
with a reflectivity of  $\Gamma_{R_o}$.   

The symmetrized field measured at the Earth's surface is then of the form
\begin{eqnarray}
\hat{p}^{\rm sym}_{l,m} (\omega,R_o)  = 
 \hat{p}_{l,m} (\omega,R_o) - \overline{\hat{p}}_{l,m} (\omega,R_o)  
  =  \hat{F}^{\rm e1}_{l,m} (\omega) (1-\Gamma_{R_o})^{2}   \hat{\cal S}_{l}(\omega,R_o)
  ,
  \label{eq:symp1}
 \end{eqnarray}
 with $\hat{F}^{\rm e1}_{l,m}$ the mode-dependent effective source trace
  and the  factor   $(1-\Gamma_{R_o})^{2}$
is an effective transmittivity factor through the  Earth's surface
(see   Proposition \ref{prop2}).

\subsubsection{Distributed Random  Source Field.}
  
 This  is our main configuration and the main result is that under ``ideal''  circumstances  we can
 recreate the surface Earth response function  in Eq.  (\ref{eq:SO}) via forming correlations. 
 Specifically, we model the source field as a random
 field  delta-correlated in space located below the surface:
 \begin{eqnarray}
\fl
\EE \big[ {f}(t,r,\theta,\varphi){f}(t',r',\theta',\varphi') ] = F (  {t-t'}  )
K(r)\delta(r-r') \sin(\theta)^{-1} \delta(\theta-\theta') \delta(\varphi-\varphi') .
\label{eq:conf}
\end{eqnarray}
 
We then have for the empirical correlation
\begin{equation*}
\label{eq:empiricalcov0}
C^{T}_{l,m}(t,R_o) = \frac{1}{T} \int_0^T p_{l,m}(t',R_o)p_{l,m}(t'+t,R_o)dt',
\end{equation*}
that for  $T$ large: 
\begin{equation}
 \hat{C}^{T}_{l,m}(\omega,R_o) =  
 \hat{F}^{\rm e2}_{l,m} (\omega) (1+\Gamma_{R_o})^{2}   \hat{\cal S}_{l}(\omega,R_o) ,
 \label{eq:CT1}
\end{equation}
with  again $\hat{F}^{\rm e2}_{l,m}$ an effective source trace (see Proposition \ref{prop:statauto}).  
  
This  generalizes  the classic daylight imaging result to the context of the spherical
Earth.  
By looking at (\ref{eq:CT1}) we can immediately deduce that the eigenfrequencies
-for which $ \hat{\cal S}_{l}(\omega,R_o) $ essentially blows up- can be extracted from the  empirical correlation  $\hat{C}^{T}_{l,m}(\omega,R_o) $.
Moreover, if one wants to reconstruct the Earth scattering operator 
$ \hat{\cal S}_{l}(\omega,R_o) $ quantitatively, then some additional assumptions have to be made on the
noise source distribution so that one can extract $ \hat{\cal S}_{l}(\omega,R_o) $ from the product 
$ F^{\rm e2}_{l,m} (\omega)  \hat{\cal S}_{l}(\omega,R_o) $.
Such an additional assumption is proposed and discussed at the end of Section \ref{sec:correlation}:
If the spatial support of the noise source is localized in a small annulus below the surface, then
the autocorrelation is directly related to the scattering operator via a classical seismic interferometry formula
and one can extract $ \hat{\cal S}_{l}(\omega,R_o) $ from $ \hat{C}^{T}_{l,m}(\omega,R_o)$
up to a multiplicative frequency-dependent function (see Proposition \ref{prop4}).

\subsection{Robustness of Main Result with Random Source Field}

The assumptions that the Earth is radial, that the noise source field is statistically homogeneous in angles 
(Eq.~(\ref{eq:conf})),  and that the wave field
is recorded  over the whole surface may seem quite restrictive. However, as we discuss
in detail in Section \ref{sec:rob}  the  analysis indicates that our results are surprisingly  robust 
with respect to these assumptions. These observations are consistent with  the
success of helioseismology where the ``source'' 
waves are generated by the turbulence in the convection zone immediately beneath the Sun's surface
\cite{sun}
and the eigenfrequencies  can be relatively robustly observed  and giving information
about the radial variation of the Sun's parameters \cite{sun1,sun2}. 

We briefly   comment on the robustness and refer to Section \ref{sec:rob} 
for a detailed discussion. 

\begin{itemize}
\item
{\it Full surface measurement aperture}:
In fact,  exploiting spherical
symmetry \redc{our result is essentially unchanged with very limited measurements},
see Section \ref{rob1}.  

\item
{\it Full surface source aperture}: 
The assumption of noise sources whose angular distribution is uniform can be slightly relaxed,
see Section \ref{rob2}.

\item
{\it A radial Earth}: We can allow for small undulations  in the velocity model beyond the radial case, see Section \ref{rob3}. 

\item
{\it Noisy measurements}:   The fact  that we average the observations over different locations
makes the scheme robust with respect to additive measurement noise, see Section
\ref{rob4}. 

\end{itemize}  
 
\section{Wave Field Asymptotics in a Stratified  Ball}

\subsection{Helmholtz Equation in the Ball}

\redc{In this subsection we formulate the Helmholtz-type equation that the mode amplitudes satisfy 
with the appropriate boundary conditions.}
In spherical coordinates $(r,\theta,\varphi)$, the time-harmonic wave field
\begin{eqnarray*}
\hat{p}(\omega,r,\theta,\varphi ) = \int p(t,r,\theta,\varphi) e^{i \omega t} dt 
\end{eqnarray*}
\redc{satisfies Eq.~(\ref{eq:main2})}
%\begin{eqnarray}\label{eq:main}
%\frac{\partial}{\partial r} r^2 \frac{\partial}{\partial r} \hat{p} +\frac{1}{\sin \theta}
%\frac{\partial}{\partial \theta} \sin \theta \frac{\partial}{\partial \theta} \hat{p}
%+\frac{1}{\sin^2 \theta} \frac{\partial^2 }{\partial \varphi^2} \hat{p}
%+ \frac{\omega^2 r^2}{c^2(r)} \hat{p} = \hat{f}(\omega,r,\theta,\varphi) ,
%\end{eqnarray}
for $(r,\theta,\varphi) \in (0,R_o) \times (0,\pi) \times(0,2\pi)$,
where $c(r)$ is the velocity model and $\hat{f}(\omega,r,\theta,\varphi)$
is the source term in the frequency domain.
\redc{If we consider the situation when this model comes from the 
acoustic wave  equation, then this means that the density is 
constant and only the bulk modulus is heterogeneous so that 
we can view $c(r)$ as the independent variable as we do here. }
We consider two types of sources in this paper:\\
- the source corresponding to seismology is point-like and located just below the surface:
\begin{eqnarray}
\label{eq:source1}
 {f}(t,r,\theta,\varphi) =  f(t) g(\theta,\varphi) \delta(r -  {R_s}) ,
\end{eqnarray}
where $ {f}(t) $ is a short pulse. We will discuss this case in detail in Section  \ref{sec:reflection}.    \\
- the source corresponding to daylight imaging is localized in an annulus and emits random noise that is  uncorrelated in space and 
stationary in time:
\begin{eqnarray}
\fl
\EE \big[ {f}(t,r,\theta,\varphi) ] = 0,\label{eq:source2a} \\ 
\fl 
\EE \big[ {f}(t,r,\theta,\varphi){f}(t',r',\theta',\varphi') ] = F(t-t')
K(r)\delta(r-r') \sin(\theta)^{-1} \delta(\theta-\theta') \delta(\varphi-\varphi'),
\label{eq:source2}
\end{eqnarray}
where $F(t-t')$ is the time correlation function of the noise sources (its Fourier transform is the power spectral density) 
and $K(r)$ is a smooth function  compactly  supported in an annulus 
$r \in [ R_s,R_o]$
below the surface.
We will discuss this case in detail in Section~\ref{sec:correlation}.  

We have Neumann (traction-free) boundary condition at the surface $r=R_o$,
which corresponds to a Robin-type boundary condition:
\begin{eqnarray}
\label{eq:tractionfree}
\big( \partial_r \hat{p} - r^{-1} \hat{p} \big)_{r=R_o} = 0 .
\end{eqnarray}
\redc{In fact, we will consider in the following a modified version (see Eq.~(\ref{eq:dissbc})) that takes 
into account boundary dissipation  \cite{lagnese83,triggiani89} so that when an interior wave hits the 
surface not all energy is reflected back into the interior.} 

At the center $r=0$ the condition is that the field should not have any singularity.

We expand the wave field  in  spherical harmonics \redc{as in Eq.~(\ref{eq:fieldmode1a0}).}
%\begin{eqnarray*}
%\hat{p}(\omega,r,\theta,\varphi) = \sum_{l=0}^\infty \sum_{m=-l}^l Y_{l,m}(\theta,\varphi) \hat{p}_{l,m} (\omega,r) .
%\end{eqnarray*}
The spherical harmonics
\begin{eqnarray}
\label{eq:Ydef}
Y_{l,m}(\theta,\varphi)= \sqrt{ \frac{(2l+1)(l-m)!}{4 \pi (l+m)!}} P^m_l(\cos \theta) e^{i m \varphi } ,
\end{eqnarray}
with $P^m_l$ being the associated Legendre polynomials,
form a complete orthonormal set in  the space $L^2(\sin \theta d\theta d\varphi)$ 
and they satisfy:
\begin{eqnarray*}
&&\int_0^\pi \int_0^{ 2\pi} \overline{Y_{l,m}(\theta,\varphi)}  Y_{l',m'}(\theta,\varphi) \sin \theta d\varphi d \theta = \delta_{ll'} \delta_{mm'} ,\\
&&\sum_{l,m} \overline{ Y_{l,m}(\theta,\varphi)}  Y_{l,m}(\theta',\varphi') \sin \theta' = 
\delta(\theta-\theta') \delta (\varphi-\varphi') .
\end{eqnarray*}
The mode amplitudes $\hat{p}_{l,m} (\omega, r) $ \redc{defined by (\ref{eq:fieldmode1a00})} satisfy the system of
second-order ODEs:
\begin{eqnarray}\label{eq:coupling}
\frac{\partial}{\partial r} r^2 \frac{\partial}{\partial r} \hat{p}_{l,m}
-l(l+1)    \hat{p}_{l,m}
+  \frac{\omega^2}{c^2(r)} r^2   \hat{p}_{l,m} = \hat{f}_{l,m}(\omega,r),
\end{eqnarray}
for $r \in (0,R_o)$,
where
\begin{eqnarray}
\label{eq:source}
 \hat{f}_{l,m}(\omega,r) &=&  \int_0^\pi \int_0^{ 2\pi} \overline{Y_{l,m}(\theta,\varphi) }
 \hat{f}(\omega,r,\theta,\varphi) \sin \theta d\varphi d\theta  .
   %  \\      &=& F_{l,m} \delta(r-R_s) . 
\end{eqnarray}
At the surface $R_o$ the mode amplitudes satisfy the Robin condition:
\begin{eqnarray}
\label{eq:robin1}
\partial_r \hat{p}_{l,m}(\omega,R_o) - R_o^{-1} \hat{p}_{l,m}(\omega,R_o)=0 .
\end{eqnarray}
At the center $r=0$ the condition is that the mode amplitude $\hat{p}_{l,m}$ should not have any singularity.
This condition can be made more explicit if we assume that the velocity is homogeneous  
in a small ball with radius $R_\delta$ at the center.
Then the mode amplitude $\hat{p}_{l,m}$ satisfies 
\begin{eqnarray*}
\frac{\partial}{\partial r} r^2 \frac{\partial}{\partial r} \hat{p}_{l,m}
-l(l+1)    \hat{p}_{l,m}
+ \frac{\omega^2}{c^2(0) } r^2   \hat{p}_{l,m} = 0,
\end{eqnarray*}
for $r \in (0,R_\delta)$. Since the only regular solution of this second-order ODE is
the spherical Bessel function, we get that $\hat{p}_{l,m}$ must be equal to
$ j_l \big(  \omega  r /c(0) \big) $,
for $r\in (0,R_\delta)$,
up to multiplicative constant.  We consider here $\omega>0$.
Therefore an explicit boundary condition is that the mode amplitude should satisfy
the Robin condition at $r=R_\delta$:
\begin{eqnarray}
\label{eq:robin2}
j_l\big(  \frac{\omega}{c(0)} R_\delta \big) \frac{\partial}{\partial r} \hat{p}_{l,m}(\omega,R_\delta) 
-  \frac{\omega}{c(0)}  j_l' \big( \frac{\omega}{c(0)}  R_\delta \big) \hat{p}_{l,m}(\omega,R_\delta) =0 .
\end{eqnarray}

\subsection{Radial Wavefield Decomposition}
\redc{In this subsection we introduce a mode decomposition that turns out to be useful in the forthcoming 
high-frequency analysis.}
Let $c_o(r)$ be the smooth component of the speed of propagation:
\begin{eqnarray*}
\frac{1}{c^2(r)} = \frac{1}{c_o^2(r)} + V(r)  ,
\end{eqnarray*}
where $V(r)$ contains the rapidly varying component responsible for scattering.
We introduce two linearly independent solutions $A_l(\omega,r)$ and $B_l(\omega,r)$  of  
the second-order ODE:
\begin{eqnarray}\label{eq:Al}
\frac{\partial}{\partial r} r^2 \frac{\partial}{\partial r} A_l
-l(l+1)   A_l
+ \frac{\omega^2}{c_o^2(r)} r^2  A_l = 0 ,
\end{eqnarray}
for $r \in (R_\delta,R_o)$.
We do not require the solutions to satisfy any boundary condition, but we require 
these independent solutions to have a Wronskian equal to one: 
\begin{eqnarray}
\label{wrons}
r^2 \frac{\partial A_l}{\partial r} B_l - r^2 A_l \frac{\partial B_l}{\partial r}  =1 ,
\end{eqnarray}
for $r \in (R_\delta,R_o)$, which is possible by
(\ref{eq:Al}) \redc{and indeed ensures linear independence}. We will determine the choice for solutions $A_l(\omega,r)$ and $B_l(\omega,r)$  later on.

We define:
\begin{eqnarray}\label{eq:an10}
 \hat{\alpha}_{l,m}(\omega,r) &=& r^2 \big( \partial_r \hat{p}_{l,m}(\omega,r) \big) B_l(\omega,r) - r^2 \hat{p}_{l,m}(\omega,r)\big( \partial_r B_l(\omega,r)\big) ,\\
 \hat{\beta}_{l,m}(\omega,r) &=& - r^2 \big( \partial_r \hat{p}_{l,m}(\omega,r) \big) A_l(\omega,r) + r^2 \hat{p}_{l,m}(\omega,r) \big( \partial_r A_l(\omega,r) \big) .\label{eq:an20}
\end{eqnarray}
The mode amplitude $\hat{p}_{l,m}$ can then be written in the form
\begin{eqnarray}\label{eq:an1}
\hat{p}_{l,m}(\omega,r) =  A_l(\omega,r)  \hat{\alpha}_{l,m}(\omega,r)+ B_l(\omega,r) \hat{\beta}_{l,m}(\omega,r) .
\end{eqnarray}

Differentiating  (\ref{eq:an10})  and   (\ref{eq:an20}) we get using (\ref{eq:Al}):
\begin{eqnarray}
\partial_r \hat{\alpha}_{l,m}  &=&
\Big( \partial_r \big( r^2 \big( \partial_r \hat{p}_{l,m} \big)\big) \Big) B_l  
  -   \hat{p}_{l,m} \Big( \partial_r \big( r^2   \big( \partial_r B_l \big)  \big)  \Big)
  \\ &=& \nonumber 
     \Big(  \partial_r \big( r^2 \big( \partial_r \hat{p}_{l,m} \big)\big) 
  -l(l+1)   \hat{p}_{l,m} 
+ \frac{\omega^2}{c_o^2(r)} r^2 \hat{p}_{l,m}   \Big) B_l   ,
  \\
\partial_r \hat{\beta}_{l,m}  &=&
 - \Big(  \partial_r \big( r^2 \big( \partial_r \hat{p}_{l,m} \big)\big) \Big) A_l  
  +   \hat{p}_{l,m}\Big( \partial_r \big( r^2   \big( \partial_r A_l \big)  \big)  \Big)
  \\ &=& \nonumber 
   -  \Big(  \partial_r \big( r^2 \big( \partial_r \hat{p}_{l,m} \big)\big) 
  - l(l+1)   \hat{p}_{l,m} 
+ \frac{\omega^2}{c_o^2(r)}  r^2  \hat{p}_{l,m}  \Big) A_l   ,
\end{eqnarray}  \redc{
and thus the amplitudes $ \hat{\alpha}_{l,m}$ and $ \hat{\beta}_{l,m}$  satisfy
\begin{eqnarray}\label{eq:an2}
A_l(\omega,r) \partial_r \hat{\alpha}_{l,m}(\omega,r) + B_l(\omega,r) \partial_r \hat{\beta}_{l,m}(\omega,r)=0  .
\end{eqnarray}}

Therefore, using (\ref{eq:coupling})
the mode amplitudes $ \hat{\alpha}_{l,m}$ and $ \hat{\beta}_{l,m}$
satisfy the system of first-order ODEs:
\begin{eqnarray}
\nonumber
\frac{\partial \hat{\alpha}_{l,m}}{\partial r} (\omega,r)&=&
- \omega^2 V(r) r^2  B_l(\omega,r) 
 \big[ A_{l}(\omega,r)  \hat{\alpha}_{l,m}(\omega,r)+ B_{l}(\omega,r) \hat{\beta}_{l,m}(\omega,r)\big] \\
\label{sys:2a} && +B_l(\omega,r)\hat{f}_{l,m} , \\
\nonumber
\frac{\partial \hat{\beta}_{l,m}}{\partial r} (\omega,r)&=&
  \omega^2 V(r) r^2 A_l(\omega,r)  
 \big[ A_{l}(\omega,r)  \hat{\alpha}_{l,m}(\omega,r)+ B_{l}(\omega,r) \hat{\beta}_{l,m}(\omega,r)\big]\\
 && -A_l(\omega,r)\hat{f}_{l,m} ,
\label{sys:2b}
\end{eqnarray}
for $r \in (R_\delta,R_o)$.

At the surface $R_o$ we find  from Eq. (\ref{eq:robin1}) that the mode amplitudes satisfy the linear relation:
\begin{eqnarray}
\nonumber
&&
 \big( \partial_r A_l(\omega,R_o) - R_o^{-1} A_l(\omega,R_o)\big) \hat{\alpha}_{l,m} (\omega,R_o)  \\
 &&
 +\big( \partial_r B_l(\omega,R_o)   - R_o^{-1}   B_l(\omega,R_o) \big) \hat{\beta}_{l,m} (\omega,R_o)=0  .
\end{eqnarray}
(note that this relation is not trivial as $\partial_r A_l - R_o^{-1} A_l$ and $\partial_r B_l   - R_o^{-1}   B_l$
cannot \redc{both be zero}  since the Wronskian of $A_l$ and $B_l$  is one).\\
At the surface $r=R_\delta$  we find from Eq. (\ref{eq:robin2}) that the mode amplitudes satisfy the linear relation:
\begin{eqnarray}
\nonumber
\Big( j_l\big(  \frac{\omega}{c(0)}  R_\delta \big)    \partial_r A_l(\omega,R_\delta)
-    \frac{\omega}{c(0)}  j_l' \big(  \frac{\omega}{c(0)} R_\delta \big) A_l(\omega,R_\delta) \Big)
\hat{\alpha}_{l,m}(\omega,R_\delta) \\
+\Big( j_l\big(  \frac{\omega}{c(0)}  R_\delta \big)    \partial_r B_l(\omega,R_\delta)
-  \frac{\omega}{c(0)}  j_l' \big(   \frac{\omega}{c(0)}  R_\delta \big) B_l(\omega,R_\delta) \Big)
 \hat{\beta}_{l,m}
(\omega,R_\delta) =0 .
\end{eqnarray}
(note that this relation is again not trivial since  the two coefficients within the big parentheses
cannot \redc{both be zero  since the Wronskian of $A_l$ and $B_l$ is one and since 
moreover  $j_l$   and $j_l'$  cannot both be zero}).

\subsection{High-frequency Asymptotics}
\label{sec:HF}
\redc{In this subsection we consider a frequency such that the radius of the Earth is much larger than the corresponding
 wavelength and we carry out high-frequency asymptotic expansions.}
%The basic scaling situation that we consider is a high-frequency situation.
More precisely, (in the assumed  non-dimensionalized coordinates)
we assume that the  surface radius $R_o$ is an order one quantity and that the 
source wavelength is small.
We therefore here consider a high frequency of the form
\begin{eqnarray*}
 \frac{\omega}{\eps} ,
\end{eqnarray*} \redc{ 
and accordingly introduce a WKB type parameterization with respect to
spherically propagating waves. }
For any $l$ we parametrize  the solutions $A_l^\eps$ and $B_l^\eps$ as 
\begin{eqnarray}
\label{eq:modeahf1a}
A_l^\eps(\omega,r) &=& \frac{\eps^{1/2}}{\sqrt{\omega}} \tilde{A}_l(r) \exp \Big( i \frac{\omega \tau(r,R_o)}{\eps} \Big) \big(1 + O(\eps)\big),\\
B_l^\eps(\omega,r) &=& \frac{\eps^{1/2}}{\sqrt{\omega}} \tilde{B}_l(r) \exp \Big( - i \frac{\omega \tau(r,R_o)}{\eps} \Big) \big(1 + O(\eps)\big),  
\label{eq:modeahf1b}
\end{eqnarray}
where $\tau(r,R_o)$ is the travel time (obtained from the eikonal equation):
\begin{eqnarray}
\label{eq:deftau}
\tau(r,R_o) = \int_r^{R_o}\frac{1}{c_o(r')} dr' ,
\end{eqnarray}
and the amplitudes $\tilde{A}_l(r) $ and $\tilde{B}_l(r) $ satisfy the frequency-independent 
transport equations:
\begin{eqnarray*}
\frac{2r^2}{c_o(r)}  \partial_r \tilde{A}_l(r)   +\partial_r \Big(\frac{r^2}{c_o(r)}\Big) \tilde{A}_l(r)  =0 ,
\end{eqnarray*}
which can be integrated as
\begin{eqnarray}
\label{eq:modeahf2}
\tilde{A}_l(r) = \tilde{A}_{l0} \frac{\sqrt{c_o(r)}}{ r} , \quad \quad  
\tilde{B}_l(r) = \tilde{B}_{l0} \frac{\sqrt{c_o(r)}}{  r} .
\end{eqnarray}
The mode $A_l^\eps$ is a down-going mode, while $B_l^\eps$ is an up-going mode.
In order to satisfy the condition (\ref{wrons}) that the Wronskian is one,  we can take:
\begin{eqnarray}
\label{def:a0}
 \tilde{A}_{l0}= \tilde{B}_{l0} = 
  \frac{e^{i \frac{\pi}{4}}}{\sqrt{2}} .
\end{eqnarray}

\begin{remark}  
The mode profiles are independent of $l$ to leading order in the regime $\eps \to 0$
as long as $l$ is of order one. 
  There are, however, corrective terms in the WKB expansion that depend on $l$.
Indeed, by computing higher order terms in  the WKB expansion, it is easy to check that
the modes can be expanded as 
\begin{eqnarray*}
\fl
A_l^\eps(\omega,r) &=&  \frac{\eps^{1/2}}{\sqrt{\omega}} \tilde{A}_l(r) \exp \Big( i \frac{\omega \tau(r,R_o)}{\eps} \Big) \\
\fl  && \times \Big( 1  - \frac{\eps}{2 i\omega} \int_r^{R_o} \frac{\sqrt{c_o(r')}}{r'}
 \Big( \partial_{r'} {r'}^2 \partial_{r'} -l(l+1) \Big) \Big(  \frac{\sqrt{c_o(r')}}{r'} \Big)   dr'
+O(\eps^{2}) 
\Big) ,\\
\fl
B_l^\eps(\omega,r) &=& \frac{\eps^{1/2}}{\sqrt{\omega}} \tilde{B}_l(r) \exp \Big( - i \frac{\omega \tau(r,R_o)}{\eps} \Big)\\
\fl && \times \Big( 1  + \frac{\eps}{2 i\omega} \int_r^{R_o} \frac{\sqrt{c_o(r')}}{r'}
\Big( \partial_{r'} {r'}^2 \partial_{r'} -l(l+1) \Big) \Big(  \frac{\sqrt{c_o(r')}}{r'} \Big)   dr'
+O(\eps^{2}) 
\Big) ,
\end{eqnarray*}
or equivalently
\begin{eqnarray*}
\fl
A_l^\eps(r) &=&  \frac{\eps^{1/2}}{\sqrt{\omega}}  \tilde{A}_l(r) \exp \Big( i \frac{\omega \tau_l^\eps(r,R_o,\omega)}{\eps} \Big) 
\Big(1 +O(\eps^{2}) 
\Big) ,\\
\fl
B_l^\eps(r) &=& \frac{\eps^{1/2}}{\sqrt{\omega}} \tilde{B}_l(r) \exp \Big( - i \frac{\omega \tau_l^\eps (r,R_o,\omega)}{\eps} \Big)
\Big(1 +O(\eps^{2}) \Big) ,
\end{eqnarray*}
with
\begin{eqnarray*}
\fl
\tau_l^\eps(r,R_o,\omega)= \tau(r,R_o) - \frac{\eps^{2}}{2 \omega^2} \int_r^{R_o} \frac{\sqrt{c_o(r')}}{r'}
 \Big( \partial_{r'} {r'}^2 \partial_{r'} -l(l+1) \Big) \Big(  \frac{\sqrt{c_o(r')}}{r'} \Big)   dr' .
\end{eqnarray*}
It is indeed possible to carry out the forthcoming analysis with these refined expressions,
but below we shall continue with the leading order phase $\tau(r,R_o)$  as this is sufficient to
characterize the leading behavior for our quantities of interest, but note here that 
\redc{the leading-order terms may take a different form} when $l={\mathcal O}(\eps^{-1/2})$. 

\end{remark}

\redc{The rapidly varying component of the speed of propagation may vary at the scale $\eps$, so we denote it by $V^\eps(r)$
so as to remember it may depend on $\eps$.}
Then, substituting the mode expressions (\ref{eq:modeahf1a}-\ref{eq:modeahf1b}) with (\ref{eq:deftau}-\ref{eq:modeahf2}-\ref{def:a0})
into the system (\ref{sys:2a}-\ref{sys:2b}),
the system of first-order ODEs for the mode amplitudes reads:
\begin{eqnarray}
\label{sys:3a}
\frac{\partial \hat{\alpha}_{l,m}^\eps}{\partial r} (\omega,r)&=&
- i \frac{\omega}{2 \eps}c_o(r)V^\eps(r)
  \big[   \hat{\alpha}_{l,m}^\eps(\omega,r)+ e^{-2i {\omega\tau(r,R_o)}/{\eps}} \hat{\beta}_{l,m}^\eps(\omega,r)\big] ,\\
\frac{\partial \hat{\beta}_{l,m}^\eps}{\partial r} (\omega,r)&=&
 i \frac{\omega}{2 \eps} c_o(r) V^\eps(r)
 \big[ e^{2i {\omega\tau(r,R_o)}/{\eps}}  \hat{\alpha}_{l,m}^\eps(\omega,r)+   \hat{\beta}_{l,m}^\eps(\omega,r)\big],
\label{sys:3b}
\end{eqnarray}
for $r \in (R_\delta,R_s)$. 
%Here we write  $V=V^\eps$ to allow
%the medium fluctuations to depend on the small parameter $\eps$,
%in particular the medium fluctuations may  be rapidly varying.

The mode amplitudes also satisfy jump and boundary conditions:\\
At the surface $R_o$ the mode amplitudes should satisfy the linear relation:
\begin{eqnarray}
\label{bc:3e0}
 \hat{\alpha}_{l,m}^\eps (\omega,R_o) - \hat{\beta}_{l,m}^\eps (\omega,R_o)=0 ,
\end{eqnarray}
which corresponds to a perfect traction-free boundary condition and a perfect reflection 
condition $\hat{\beta}_{l,m}^\eps (\omega,R_o)/  \hat{\alpha}_{l,m}^\eps (\omega, R_o)=1$ at the surface. 
If  we assume that there is some attenuation or dissipation at the surface, 
then the reflection condition is not equal to one exactly, 
but to some number $\Gamma_{R_o} \in (-1,1)$, and   the mode amplitudes  satisfy the linear relation:
\begin{eqnarray}
\label{bc:3e}
 \hat{\alpha}_{l,m}^\eps (\omega,R_o) - \Gamma_{R_o}  \hat{\beta}_{l,m}^\eps (\omega,R_o)=0 .
\end{eqnarray}
More exactly, if we consider the dissipation boundary condition \cite{lagnese83,triggiani89}:
\begin{eqnarray}
\label{eq:dissbc}
\big( \partial_r {p} - r^{-1} {p} \big)_{r=R_o} = - \big( \kappa \partial_t p \big)_{r=R_o},
\end{eqnarray}
instead of the traction free boundary condition (\ref{eq:tractionfree}),
then we get the boundary condition (\ref{bc:3e}) with 
\begin{eqnarray}
\label{eq:Gdef}
\Gamma_{R_o} = (1-\kappa c_o(R_o))/(1+\kappa c_o(R_o)).
\end{eqnarray}
We \redc{will} look at $\kappa=0$ or $\Gamma_{R_o}=1$ as a limiting case in the following.\\
At the surface $r=R_\delta$  the mode amplitudes satisfy the linear relation: 
\begin{eqnarray*}
\hat{\alpha}_{l,m}^\eps(\omega,R_\delta) - \hat{\beta}_{l,m}^\eps (\omega,R_\delta) e^{ -2i {\omega \tau(0,R_o)}/{\eps}} =0 ,
\end{eqnarray*}  
where we have used the asymptotic formulas for the spherical Bessel functions:
\begin{eqnarray*}
\fl
j_l(x) \simeq \frac{1}{x} \cos \Big( x - \frac{(l+1)\pi}{2} \Big), \quad \quad
j_l'(x) \simeq -\frac{1}{x} \sin \Big( x - \frac{(l+1)\pi}{2} \Big), \quad \quad
 x \gg 1 .
\end{eqnarray*}
Note that this condition
is independent on the value of $R_\delta$ (which can be expected as the 
mode amplitudes $\hat{\alpha}_{l,m}^\eps$ and $\hat{\beta}_{l,m}^\eps$
are constant within the homogeneous central small ball)
and that we can take the limit $R_\delta \to 0$:     
\begin{eqnarray}
\label{bc:3f}
 \hat{\alpha}_{l,m}^\eps(\omega,0) -  \hat{\beta}_{l,m}^\eps (\omega,0) \exp\Big( -2i \frac{\omega \tau(0,R_o)}{\eps} \Big)   =0 .
\end{eqnarray}

\begin{remark}
When there is no attenuation, i.e. when $\kappa=0$ and  $\Gamma_{R_o}=1$, then 
there is existence and uniqueness of the solution provided $\omega$ is not an eigenfrequency.
A frequency $\omega$ is an eigenfrequency if the 
linear system  (\ref{sys:3a}-\ref{sys:3b}) for $r \in (0,R_o)$
admits a non-zero solution that satisfies the two boundary conditions (\ref{bc:3e}-\ref{bc:3f})
when there is no source.
When there is attenuation then there is existence and uniqueness of the solution for any $\omega$,
\redc{this follows from the form of the propagator associated with (\ref{sys:3a}-\ref{sys:3b}), see  
Chap. 7 in  \cite{book}.}
\end{remark}

 Across a source at depth $R_s$ whose Fourier component at frequency $\omega/\eps$ has the form
\begin{eqnarray}
\label{eq:sources}
  \hat{f}_{l,m}^\eps(\omega,r)= F_{l,m}(\omega) \delta(r-R_s),
\end{eqnarray}
the mode amplitudes satisfy the jump source conditions:
\begin{eqnarray}
\nonumber
 \begin{pmatrix}
 \hat{\alpha}_{l,m}^\eps (\omega,R_s^+)\\
  \hat{\beta}_{l,m}^\eps (\omega,R_s^+) 
  \end{pmatrix}
  -  
  \begin{pmatrix}
  \hat{\alpha}_{l,m}^\eps (\omega,R_s^-)\\
     \hat{\beta}_{l,m}^\eps (\omega,R_s^-) 
    \end{pmatrix}   \\
= \eps^{1/2} F_{l,m} (\omega)  \frac{e^{i \pi /4} \sqrt{c_o(R_s)}}{R_s \sqrt{2 \omega}}
  \begin{pmatrix}  \exp \Big( - i \frac{\omega \tau(R_s,R_o)}{\eps} \Big)  \\
 -  \exp \Big( i \frac{\omega \tau(R_s,R_o)}{\eps} \Big) 
  \end{pmatrix} . 
  \label{bc:3c}
\end{eqnarray}

It is convenient to  introduce the propagator of the system (\ref{sys:3a}-\ref{sys:3b}), 
that is, the $2\times 2$ matrix ${\bf P}^\eps_l(\omega,r',r)$
solution of
\begin{eqnarray}
\label{eq:Pdef}
\fl
\frac{\partial }{\partial r} {\bf P}^\eps_l(\omega,r',r)=
 i \frac{\omega}{2 \eps} c_o(r) V^\eps (r)
 \begin{pmatrix}
 -1 & -e^{-2i {\omega\tau(r,R_o)}/{\eps}} \\
   e^{2i {\omega\tau(r,R_o)}/{\eps}}  &1
 \end{pmatrix}
 {\bf P}^\eps_l(\omega,r',r) ,
\end{eqnarray}
with $  {\bf P}^\eps_l(\omega,r',r=r') = {\bf I}$.
The matrix ${\bf P}^\eps_l(\omega,r',r)$
has the symplectic form
\begin{eqnarray}
 {\bf P}^\eps_l(\omega,r',r)= \begin{pmatrix}
  \hat{a}^\eps_l(\omega,r',r) &\overline{ \hat{b}^\eps_l(\omega,r',r) } \\
  \hat{b}^\eps_l(\omega,r',r) &\overline{ \hat{a}^\eps_l(\omega,r',r) }
  \end{pmatrix} ,  \label{eq:symplectic}
\end{eqnarray}
where $(\hat{a}^\eps_l(\omega,r',r) , \hat{b}^\eps_l(\omega,r',r) )^T$ satisfies
\begin{eqnarray}\label{eq:adyn}
\fl
\frac{\partial }{\partial r} \begin{pmatrix} \hat{a}^\eps_l(\omega,r',r) \\ \hat{b}^\eps_l(\omega,r',r)  \end{pmatrix} =
 i \frac{\omega}{2 \eps} c_o(r) V^\eps(r)
 \begin{pmatrix}
 -1 & - e^{-2i {\omega\tau(r,R_o)}/{\eps}} \\
  e^{2i {\omega\tau(r,R_o)}/{\eps}}  &1
 \end{pmatrix}
\begin{pmatrix} \hat{a}^\eps_l(\omega,r',r) \\ \hat{b}^\eps_l(\omega,r',r)  \end{pmatrix}  ,
\end{eqnarray}
starting from $  \hat{a}^\eps_l(\omega,r',r=r')=1$, $\hat{b}^\eps_l(\omega,r',r=r')=0$.
The solution satisfies the energy conservation relation 
\begin{eqnarray}
\label{eq:ener}
|\hat{a}^\eps_l(\omega,r',r)|^2-|\hat{b}^\eps_l(\omega,r',r)|^2=1 .
\end{eqnarray}

\subsection{Wave Decomposition in the Radial Earth}
\redc{In this subsection we summarize the basic high-frequency wave field decomposition in the stratified sphere.}
We define the scaled Fourier transform $\hat{f}^\eps(\omega) $ of a function $f(t)$ as
\begin{eqnarray}
\label{def:fourier}
\hat{f}^\eps(\omega) =  \frac{1}{\eps}  \int f(t) e^{i \omega t /\eps} {dt}  .
\end{eqnarray}
In spherical coordinates $(r,\theta,\varphi)$, the time-harmonic wave field
$\hat{p}^\eps(\omega,r,\theta,\varphi )  $
satisfies
\begin{eqnarray}\label{eq:main2b}
\fl
\frac{\partial}{\partial r} r^2 \frac{\partial}{\partial r} \hat{p}^\eps +\frac{1}{\sin \theta}
\frac{\partial}{\partial \theta} \sin \theta \frac{\partial}{\partial \theta} \hat{p}^\eps
+\frac{1}{\sin^2 \theta} \frac{\partial^2 }{\partial \varphi^2} \hat{p}^\eps
+ \frac{\omega^2 r^2}{\eps^{2} c^2(r)} \hat{p}^\eps =   \hat{f}^\eps(\omega,r,\theta,\varphi) ,
\end{eqnarray}
for $(r,\theta,\varphi) \in (0,R_o) \times (0,\pi) \times(0,2\pi)$,
where $c(r)$ is the velocity model and $\hat{f}^\eps(\omega,r,\theta,\varphi)$
is the source term in the frequency domain.
We assume here that 
$\hat{f}^\eps$ is supported away from the origin in $\omega$. 
The radially dependent speed parameter $c(r)$ is modeled by
\begin{eqnarray}
\frac{1}{c^2(r)} = \frac{1}{c_o^2(r)} + V^{\eps}(r)  .
\end{eqnarray}
The wave field is then decomposed as
\begin{eqnarray}
\label{eq:fieldmode1a}
\hat{p}^\eps (\omega,r,\theta,\varphi) =& \sum_{l,m} Y_{l,m}(\theta,\varphi)\hat{p}_{l,m}^\eps (\omega,r),\\
\nonumber
\hat{p}_{l,m}^\eps (\omega,r) =& {\eps^{1/2}} \frac{\sqrt{c_o(r)} e^{i\pi/4}}{r \sqrt{2 \omega}} 
\Big(\hat{\alpha}_{l,m}^\eps(\omega,r) \exp \big( \frac{i\omega \tau(r,R_o)}{\eps}\big)  \\
& +\hat{\beta}_{l,m}^\eps(\omega,r) \exp \big( - \frac{i\omega \tau(r,R_o)}{\eps}\big)  \Big)  ,
\label{eq:fieldmode1b}
\end{eqnarray}
with $Y_{l,m}$ defined in Eq.  (\ref{eq:Ydef}).

 For $0 < r' <  r  < R_o$ so that  there is no source in the interval
 $(r',r)$,  the mode amplitudes are related  
 by the propagator as
\begin{eqnarray}
  \begin{pmatrix}
  \hat{\alpha}^\eps_{l,m}(\omega,r)  \\
  \hat{\beta}^\eps_{l,m}(\omega,r)  
  \end{pmatrix} 
 =
 {\bf P}^\eps_l(\omega,r',r) 
  \begin{pmatrix}
  \hat{\alpha}^\eps_{l,m}(\omega,r')  \\
  \hat{\beta}^\eps_{l,m}(\omega,r')  
  \end{pmatrix} , \label{eq:s2}
\end{eqnarray}
 with the propagator ${\bf P}^\eps_l$ defined by   (\ref{eq:Pdef}). 
  
 At the surface $r=R_o$ the mode amplitudes  satisfy the boundary condition~(\ref{bc:3e}).  
  
At the center $r=0$ the mode amplitudes satisfy the boundary condition~(\ref{bc:3f}).

Assume a source at depth $r=R_s$  as
\begin{eqnarray}
\label{eq:source3}
 \hat{f}_{l,m}^\eps(\omega,r) 
 &=&  \int_0^\pi \int_0^{ 2\pi} \overline{Y_{l,m}(\theta,\varphi) }
 \hat{f}^\eps(\omega,r,\theta,\varphi) \sin \theta d\varphi d\theta   \\
%  &=&  \int_0^\pi \int_0^{ 2\pi}  \int_{-\infty}^\infty   \overline{Y_{l,m}(\theta,\varphi) }
%  {f}(t,r,\theta,\varphi)   e^{it\omega} dt \sin \theta d\varphi d\theta   \\
      &=& F_{l,m} (\omega) \delta(r-R_s)  ,  \label{eq:source3b}
\end{eqnarray}
then the mode amplitudes satisfy the jump conditions in Eq. (\ref{bc:3c}).

\section{Scattering in Seismology}\label{sec:reflection}
We consider a source of the form (\ref{eq:source1}) just below
the surface $R_s=R_o^{-}$
whose emission has a typical frequency of order $\eps^{-1}$:
\begin{eqnarray*}
{f}(t,r,\theta,\varphi) =  f\Big(\frac{t}{\eps}\Big) g(\theta,\varphi) \delta(r -  {R_s}) ,
\end{eqnarray*}
where $f$ is a function whose Fourier transform is supported away from the origin.
Therefore the source term defined as in (\ref{def:fourier})-(\ref{eq:source3}) reads
\begin{eqnarray*}
\hat{f}^\eps_{l,m}(\omega,r) = F_{l,m}(\omega) \delta(r -  {R_s}) ,
\end{eqnarray*}
with
\begin{eqnarray*}
 F_{l,m}(\omega)=  \hat{f} ( \omega ) {g}_{l,m},\quad \quad 
g_{l,m} = \int_0^\pi \int_0^{ 2\pi} \overline{Y_{l,m}(\theta,\varphi)} g(\theta,\varphi) \sin \theta d\theta d\varphi  .
\end{eqnarray*}
The mode amplitudes satisfy the jump source conditions:
\begin{eqnarray}
\label{eq:combiab11a}
&&\big[    \hat{\alpha}_{l,m}^\eps \big]_{R_o^-}^{R_o} A_l(\omega,R_o)
+ \big[ \hat{\beta}_{l,m}^\eps \big]_{R_o^-}^{R_o} B_l(\omega,R_o) = 0 , \\
\label{eq:combiab11b}
&&\big[    \hat{\alpha}_{l,m}^\eps \big]_{R_o^-}^{R_o} \partial_r A_l(\omega,R_o) + \big[ \hat{\beta}_{l,m}^\eps
\big]_{R_o^-}^{R_o}  \partial_r B_l(\omega,R_o) = \frac{ F_{l,m}(\omega) }{R_o^2} ,
\end{eqnarray}
\redc{We sum Eq.~(\ref{eq:combiab11a}) multiplied by $-R_o^2\partial_r B_l(\omega,R_o)$
and Eq.~(\ref{eq:combiab11b}) multiplied by $R_o^2B_l(\omega,R_o)$, and we make use of the Wronskian relation~(\ref{wrons})
to get}
\begin{eqnarray}
\big[    \hat{\alpha}_{l,m}^\eps \big]_{R_o^-}^{R_o} =  F_{l,m} (\omega) B_l(\omega,R_o) ,
\end{eqnarray}
\redc{Similarly, we sum Eq.~(\ref{eq:combiab11a}) multiplied by $R_o^2\partial_r A_l(\omega,R_o)$
and Eq.~(\ref{eq:combiab11b}) multiplied by $-R_o^2A_l(\omega,R_o)$ and we obtain}
\begin{eqnarray}
\big[ \hat{\beta}_{l,m}^\eps \big]_{R_o^-}^{R_o}   = - F_{l,m} (\omega)A_l(\omega,R_o) .
\end{eqnarray}

At the surface we ``measure'' the spherical harmonics  of  the field:
\begin{eqnarray*}
p_{l,m} (t,R_o) = \frac{1}{2\pi} \int \hat{p}_{l,m}^\eps (\omega,R_o) e^{-i  {\omega} t /\eps } d\omega  .
\end{eqnarray*}

\begin{proposition}\label{prop1}
The measured field is of the form
\begin{eqnarray}
\label{eq:eq:fieldmode1bc}
\hat{p}_{l,m}^\eps (\omega,R_o) 
&=&- \frac{i \eps c_o(R_o)}{\omega R_o^2 }
%\hat{f}(\omega) g_{l,m} 
 F_{l,m} (\omega) 
\frac{1+\Gamma_{R_o}}{2}
\frac{1 + \hat{\cal R}_{l}^\eps(\omega,R_o)}
{1- \Gamma_{R_o} \hat{\cal R}_{l}^\eps(\omega,R_o)}  ,
 \end{eqnarray}
where we have defined  \redc{the modulus one function}
\begin{eqnarray}
\label{def:Reps}
 \hat{\cal R}_{l}^\eps(\omega,R_o) 
 &=&
\frac{
\hat{b}_{l}^\eps(\omega,0,R_o)  +\overline{\hat{a}_{l}^\eps(\omega,0,R_o)}
e^{2 i \omega \tau(0,R_o)/\eps} 
}{
\hat{a}_{l}^\eps(\omega,0,R_o)  + \overline{\hat{b}_{l}^\eps(\omega,0,R_o)} 
e^{2 i \omega \tau(0,R_o)/\eps} 
}
 .
\end{eqnarray}
\end{proposition}

\debproof
We have from \redc{Eq. (\ref{eq:fieldmode1b})}:
\begin{eqnarray}
\label{eq:eq:fieldmode1bb}
\hat{p}_{l,m}^\eps (\omega,R_o) = \eps^{1/2} \frac{\sqrt{c_o(R_o)}e^{i \pi/4} }{R_o \sqrt{2\omega}} 
\big(\hat{\alpha}_{l,m}^\eps(\omega,R_o)  
+\hat{\beta}_{l,m}^\eps(\omega,R_o)   \big) .
\end{eqnarray}
\redc{The mode amplitudes $\hat{\alpha}_{l,m}^\eps$ and $\hat{\beta}_{l,m}^\eps$ satisfy
the boundary condition (\ref{bc:3e}) at $r=R_o$, 
the source jump condition (\ref{bc:3c}) at $r=R_o^-$,
the propagation equation (\ref{eq:s2}) from $r=0$ to $r=R_o^-$ (using the form  (\ref{eq:symplectic}) of the propagator matrix),
and the boundary condition (\ref{bc:3e}) at $r=0$:}
\begin{eqnarray*}
&& \hat{\alpha}_{l,m}^\eps (\omega,R_o) - \Gamma_{R_o}  \hat{\beta}_{l,m}^\eps (\omega,R_o)=0 ,\\
&& \begin{pmatrix}
 \hat{\alpha}_{l,m}^\eps (\omega,R_o)\\
  \hat{\beta}_{l,m}^\eps (\omega,R_o) 
  \end{pmatrix}
  =
  \begin{pmatrix}
  \hat{\alpha}_{l,m}^\eps (\omega,R_o^-)\\
     \hat{\beta}_{l,m}^\eps (\omega,R_o^-) 
    \end{pmatrix}   
  + \eps^{1/2}  F_{l,m}(\omega) \frac{e^{i \pi /4} \sqrt{c_o(R_o)}}{R_o \sqrt{2 \omega}}
  \begin{pmatrix} 1 \\
 - 1
  \end{pmatrix} ,\\
  &&
   \begin{pmatrix}
  \hat{\alpha}_{l,m}^\eps (\omega,R_o^-)\\
     \hat{\beta}_{l,m}^\eps (\omega,R_o^-) 
    \end{pmatrix}   
  =
  \begin{pmatrix}
  \hat{a}^\eps_{l}(\omega,0,R_o) &\overline{ \hat{b}^\eps_{l}(\omega,0,R_o)  } \\
  \hat{b}^\eps_{l}(\omega,0,R_o)  &\overline{ \hat{a}^\eps_{l}(\omega,0,R_o)  }
  \end{pmatrix}
   \begin{pmatrix}
  \hat{\alpha}_{l,m}^\eps (\omega,0)\\
     \hat{\beta}_{l,m}^\eps (\omega,0) 
    \end{pmatrix}    ,\\
    &&
     \hat{\alpha}_{l,m}^\eps(\omega,0) -  \hat{\beta}_{l,m}^\eps (\omega,0) \exp\Big( -2i \frac{\omega \tau(0,R_o)}{\eps} \Big)   =0 .
  \end{eqnarray*}
  \redc{
By solving this $6\times 6$ linear system for 
 $$
 (\hat{\alpha}_{l,m}^\eps (\omega,R_o) ,  \hat{\beta}_{l,m}^\eps (\omega,R_o),
   \hat{\alpha}_{l,m}^\eps (\omega,R_o^-),      \hat{\beta}_{l,m}^\eps (\omega,R_o^-) ,  \hat{\alpha}_{l,m}^\eps (\omega,0) , \hat{\beta}_{l,m}^\eps (\omega,0) )
   $$ we obtain the expression of $(\hat{\alpha}_{l,m}^\eps (\omega,R_o) ,  \hat{\beta}_{l,m}^\eps (\omega,R_o))$ 
   that we substitute into (\ref{eq:eq:fieldmode1bb}), which gives the desired result (\ref{eq:eq:fieldmode1bc}). }
\finproof
 
If we assume that the pulse is an even function and we consider the symmetrized response function
\begin{eqnarray}
\label{eq:responsesym}
p^{{\rm sym}}_{l,m}(t,R_o) = p_{l,m}(t,R_o)-p_{l,m}(-t,R_o) ,
\end{eqnarray}
then its Fourier transform (\ref{def:fourier}) is of the form
\begin{eqnarray*}
\hat{p}^{\eps,{\rm sym}}_{l,m}(\omega,R_o) = - \frac{i \eps c_o(R_o)}{R_o^2 \omega}
%\hat{f}(\omega) g_{l,m} 
 F_{l,m} (\omega) 
{\rm Re}\Big\{ (1+\Gamma_{R_o})
\frac{1 + \hat{\cal R}_{l}^\eps(\omega,R_o)}
{1- \Gamma_{R_o} \hat{\cal R}_{l}^\eps(\omega,R_o)} \Big\} ,
\end{eqnarray*}
which is equal to
\begin{eqnarray*}
\hat{p}^{\eps,{\rm sym}}_{l,m}(\omega,R_o) = - 
 \frac{i \eps c_o(R_o)}{ R_o^2 \omega}
%\hat{f}(\omega) g_{l,m} 
 F_{l,m} (\omega) 
\frac{\big(1 -|\Gamma_{R_o}|^2\big) 
{\rm Re}\big\{ 1 + \hat{\cal R}_{l}^\eps(\omega,R_o) \big\}
}
{\big| 1- \Gamma_{R_o} \hat{\cal R}_{l}^\eps(\omega,R_o) \big|^2} .
\end{eqnarray*}
Using $|1+\hat{\cal R}_{l}^\eps(\omega,R_o) |^2=2 {\rm Re}\big\{ 1 + \hat{\cal R}_{l}^\eps(\omega,R_o) \big\}$,
we finally establish the main result of this section.
\begin{proposition}\label{prop2}
The symmetrized response function is of the form:
\begin{eqnarray}
\label{eq:symrep2}
\hat{p}^{\eps,{\rm sym}}_{l,m}(\omega,R_o) = - \frac{i \eps c_o(R_o)}{2 R_o^2 \omega}
%\hat{f}(\omega) g_{l,m} 
 F_{l,m} (\omega) \big(1 -|\Gamma_{R_o}|^2\big) 
\frac{
\big| 1 + \hat{\cal R}_{l}^\eps(\omega,R_o) \big|^2
}
{\big| 1- \Gamma_{R_o} \hat{\cal R}_{l}^\eps(\omega,R_o) \big|^2} ,
\end{eqnarray}
where $\hat{\cal R}_{l}^\eps(\omega,R_o) $ is defined by (\ref{def:Reps}).
\end{proposition}

The scattering operator 
\begin{eqnarray}
\label{eq:refop}
\frac{\big| 1 + \hat{\cal R}_{l}^\eps(\omega,R_o) \big|^2}{|1-\Gamma_{R_o} \hat{\cal R}_{l}^\eps(\omega,R_o) |^2}
\end{eqnarray}
 is the quantity of interest that characterizes 
the medium and that can be extracted from standard seismology.
Recall  that $ \hat{\cal R}_{l}^\eps(\omega,R_o) $ has modulus one.
Frequencies such that $\hat{\cal R}_{l}^\eps(\omega,R_o)=1$ 
correspond to the eigenmodes, and we remark that then
\begin{eqnarray*}
\big(1 -|\Gamma_{R_o}|^2 \big)
\frac{ \big| 1 + \hat{\cal R}_{l}^\eps(\omega,R_o) \big|^2}
{\big| 1- \Gamma_{R_o} \hat{\cal R}_{l}^\eps(\omega,R_o) \big|^2}  \simeq \frac{4 }{\kappa c_o(R_o)} ,
\end{eqnarray*}
as $\kappa\to 0$ 
in view of the relation  (\ref{eq:Gdef}).

\section{Daylight Imaging}
\label{sec:correlation}%
We consider a noise source  as in  
(\ref{eq:source2a})  and  (\ref{eq:source2}) with  a correlation function of the form
\begin{eqnarray}\label{eq:ds}
\fl
\EE \big[ {f}(t,r,\theta,\varphi){f}(t',r',\theta',\varphi') ] = F\Big( \frac{t-t'}{\eps}\Big)
K(r)\delta(r-r') \sin(\theta)^{-1} \delta(\theta-\theta') \delta(\varphi-\varphi'),
\end{eqnarray}
which means that the power spectrum of the noise source contains frequencies of the order of $\eps^{-1}$.
Under these conditions the source term $\hat{f}_{l,m}^\eps(\omega,r)$ defined by 
(\ref{def:fourier}) and (\ref{eq:source3}) is a random process with mean zero and covariance
function  
\begin{eqnarray}\label{eq:spec}
\EE \big[ \hat{f}_{l,m}^\eps(\omega,r) \overline{\hat{f}_{l',m'}^\eps(\omega',r')} \big] = 2\pi
\hat{F}(\omega) K(r) \delta(r-r') \delta(\omega-\omega')   \delta_{ll'} \delta_{mm'}  .
\end{eqnarray}
The field at the surface is $p(t,R_o,\theta,\varphi)$ and we compute its empirical
autocorrelation components associated with the spherical harmonics:
\redc{
\begin{equation}
\label{eq:empiricalcov}
C^{T}_{l,m}(t,R_o) = \frac{1}{T} \int_0^T p_{l,m}(t',R_o)p_{l,m}(t'+t,R_o)dt',
\end{equation}
for
\begin{eqnarray}
 \label{eq:empiricalcov2}
p_{l,m}(t,R_o) = \int_0^\pi \int_0^{ 2\pi} p(t,R_o,\theta,\varphi) \overline{ Y_{l,m}(\theta,\varphi) }
\sin \theta d\varphi d \theta  .
\end{eqnarray} 
In practice, the signals $p_{l,m}(t,R_o)$ are computed from the signals $(p(t, R_o, \theta_j,\varphi_j))_{j=1,\ldots,N}$,
recorded by a collection of $N$ receivers located at $((R_o,\theta_j,\varphi_j))_{j=1,\ldots,N}$ at the surface of the Earth  via a quadrature formula.
We comment in more detail on robustness with respect to sampling in Section \ref{rob1}.
}

As $T \to \infty$ the correlations  converges in probability towards the statistical autocorrelation defined by
\begin{eqnarray}
\label{def:corstat}
C_{l,m}(t,R_o) = \EE \big[ p_{l,m}(0,R_o)p_{l,m}(t,R_o) \big] .
\end{eqnarray}
In the following proposition we study its scaled Fourier transform:
\begin{eqnarray*}
%C_{l,m}(t,R_o) = \frac{1}{2\pi} \int \hat{C}^\eps_{l,m}(\omega,R_o) 
%e^{-i  \omega  t /\eps } d\omega , \quad \quad
\hat{C}^\eps_{l,m}(\omega,R_o) =\frac{1}{\eps} \int C_{l,m}(t,R_o) e^{-i  \omega  t / \eps } dt.
\end{eqnarray*}
\begin{proposition}
\label{prop:statauto}
The statistical autocorrelation is of the form
\begin{eqnarray}
 \nonumber
&& 
\hat{C}_{l,m}^\eps (\omega,R_o) =
 \eps^2 
 \hat{F}(\omega) \frac{c_o(R_o)}{4 R_o^2   \omega^2}   
\frac{|1+\Gamma_{R_o}|^2}
{\big| 1- \Gamma_{R_o} \hat{\cal R}_{l}^\eps(\omega,R_o)\big|^2} 
\\ && \quad  \times
  \Big[ \int_0^{R_o}   \frac{c_o(R_s)}{R_s^2} K(R_s)  
\big|
\hat{\cal S}_{l}^\eps(\omega,R_s) 
+ \hat{\cal T}_{l}^\eps(\omega,R_s)
e^{-2 i \omega \tau(R_s,R_o)/\eps} 
\big|^2 dR_s \Big]   ,
 \label{eq:cor2a}
 \end{eqnarray}
with
\begin{eqnarray}
\label{eq:hatSeps}
\hat{\cal S}^\eps_{l}(\omega,R_s) &=& 
\frac{
\hat{a}_{l}^\eps(\omega,0,R_s)  + \overline{\hat{b}_{l}^\eps(\omega,0,R_s)} 
e^{2 i \omega \tau(0,R_o)/\eps} 
}
{
\hat{a}_{l}^\eps(\omega,0,R_o)  + \overline{\hat{b}_{l}^\eps(\omega,0,R_o)} 
e^{2 i \omega \tau(0,R_o)/\eps} 
} ,
\\
\label{eq:hatTeps}
\hat{\cal T}^\eps_{l}(\omega,R_s) &=& 
\frac{
\hat{b}_{l}^\eps(\omega,0,R_s)
 + \overline{\hat{a}_{l}^\eps(\omega,0,R_s)} 
e^{2 i \omega \tau(0,R_o)/\eps} 
}
{
\hat{a}_{l}^\eps(\omega,0,R_o)  + \overline{\hat{b}_{l}^\eps(\omega,0,R_o)} 
e^{2 i \omega \tau(0,R_o)/\eps} 
} .
\end{eqnarray}
\end{proposition}

\debproof 
We have
\begin{eqnarray*}
\fl
\hat{p}_{l,m}^\eps (\omega,R_o) =
\int_0^{R_o} \eps^{1/2} \frac{\sqrt{c_o(R_o)}e^{i \pi/4} }{R_o \sqrt{2\omega}} 
  \hat{f}_{l,m}^\eps(\omega,R_s)
\big(\hat{\alpha}_{l,m}^\eps(\omega,R_o; R_s)  
+\hat{\beta}_{l,m}^\eps(\omega,R_o; R_s)   \big) d R_s ,
\end{eqnarray*}
where $\hat{\alpha}_{l,m}^\eps(\omega,r;R_s) $ and $\hat{\beta}_{l,m}^\eps(\omega,r;R_s) $ are the mode amplitudes corresponding 
to a unit-amplitude spherical harmonic point source at $r=R_s$.
We get then in view of~(\ref{eq:spec}) and~(\ref{def:corstat}): 
\begin{eqnarray*}
\fl
\hat{C}_{l,m}^\eps (\omega,R_o) 
&=& \EE\Big [ \hat{p}_{l,m}^\eps (\omega,R_o) \frac{1}{2\pi} \int \overline{\hat{p}_{l,m}^\eps (\omega',R_o)} d\omega'\Big] \\
\fl &=& \eps \hat{F}(\omega)  \frac{c_o(R_o)}{2 R_o^2 |\omega|} 
 \int_0^{R_o} K(R_s) \big|\hat{\alpha}_{l,m}^\eps(\omega,R_o;R_s) +\hat{\beta}_{l,m}^\eps(\omega,R_o;R_s) \big|^2  dR_s .
\end{eqnarray*}
The mode amplitudes $\hat{\alpha}_{l,m}^\eps(\omega,r;R_s) $ and $\hat{\beta}_{l,m}^\eps(\omega,r;R_s) $ satisfy
the boundary, propagation, and source jump conditions:
\begin{eqnarray*}
\fl \hat{\alpha}_{l,m}^\eps (\omega,R_o;R_s) - \Gamma_{R_o}  \hat{\beta}_{l,m}^\eps (\omega,R_o;R_s)=0 ,\\
\fl
   \begin{pmatrix}
  \hat{\alpha}_{l,m}^\eps (\omega,R_o;R_s)\\
     \hat{\beta}_{l,m}^\eps (\omega,R_o;R_s) 
    \end{pmatrix}   
  =
  \begin{pmatrix}
  \hat{a}^\eps_{l}(\omega,R_s,R_o) &\overline{ \hat{b}^\eps_{l}(\omega,R_s,R_o)  } \\
  \hat{b}^\eps_{l}(\omega,R_s,R_o)  &\overline{ \hat{a}^\eps_{l}(\omega,R_s,R_o)  }
  \end{pmatrix}
   \begin{pmatrix}
  \hat{\alpha}_{l,m}^\eps (\omega,R_s^+;R_s)\\
     \hat{\beta}_{l,m}^\eps (\omega,R_s^+;R_s) 
    \end{pmatrix}    ,\\
\fl
 \begin{pmatrix}
 \hat{\alpha}_{l,m}^\eps (\omega,R_s^+;R_s)\\
  \hat{\beta}_{l,m}^\eps (\omega,R_s^+;R_s) 
  \end{pmatrix}
  =
  \begin{pmatrix}
  \hat{\alpha}_{l,m}^\eps (\omega,R_s^-;R_s)\\
     \hat{\beta}_{l,m}^\eps (\omega,R_s^-;R_s) 
    \end{pmatrix}   
  + \eps^{1/2} \frac{e^{i \pi /4} \sqrt{c_o(R_s)}}{R_s \sqrt{2 \omega}}
  \begin{pmatrix} \exp\big( - i \frac{\omega \tau(R_s,R_o)}{\eps}\big) \\
 - \exp\big( i \frac{\omega \tau(R_s,R_o)}{\eps}\big)  
  \end{pmatrix}  ,\\
\fl 
   \begin{pmatrix}
  \hat{\alpha}_{l,m}^\eps (\omega,R_s^-;R_s)\\
     \hat{\beta}_{l,m}^\eps (\omega,R_s^-;R_s) 
    \end{pmatrix}   
  =
  \begin{pmatrix}
  \hat{a}^\eps_{l}(\omega,0,R_s) &\overline{ \hat{b}^\eps_{l}(\omega,0,R_s)  } \\
  \hat{b}^\eps_{l}(\omega,0,R_s)  &\overline{ \hat{a}^\eps_{l}(\omega,0,R_s)  }
  \end{pmatrix}
   \begin{pmatrix}
  \hat{\alpha}_{l,m}^\eps (\omega,0;R_s)\\
     \hat{\beta}_{l,m}^\eps (\omega,0;R_s) 
    \end{pmatrix}    ,\\
\fl
     \hat{\alpha}_{l,m}^\eps(\omega,0;R_s) -  \hat{\beta}_{l,m}^\eps (\omega,0;R_s) \exp\Big( -2i \frac{\omega \tau(0,R_o)}{\eps} \Big)   =0 .
  \end{eqnarray*}
By solving this linear system we get the desired result. To get the desired expression (\ref{eq:cor2a}),
we use the energy conservation relation (\ref{eq:ener}):
\begin{eqnarray*}
|\hat{a}^\eps_l(\omega,R_s,R_o)|^2-|\hat{b}^\eps_l(\omega,R_s,R_o)|^2=1 ,
\end{eqnarray*}
and we also make use of the propagation relation 
${\bf P}^\eps_l(\omega,0,R_o)={\bf P}^\eps_l(\omega,R_s,R_o) {\bf P}^\eps_l(\omega,0,R_s)$
which gives \redc{
\begin{eqnarray*}
\hat{a}^\eps_l(\omega,0,R_o) &=&
\hat{a}^\eps_l(\omega,R_s,R_o) \hat{a}^\eps_l(\omega,0,R_s) 
+\overline{\hat{b}^\eps_l(\omega,R_s,R_o) } \hat{b}^\eps_l(\omega,0,R_s)  ,\\
\hat{b}^\eps_l(\omega,0,R_o) &=&
\hat{b}^\eps_l(\omega,R_s,R_o) \hat{a}^\eps_l(\omega,0,R_s) 
+\overline{\hat{a}^\eps_l(\omega,R_s,R_o) } \hat{b}^\eps_l(\omega,0,R_s)  .
\end{eqnarray*}}
\finproof

Although the autocorrelation function (\ref{eq:cor2a}) is related to the scattering operator of interest (\ref{eq:refop}), 
it also  depends on the source distribution in a non-trivial way,
which may complicate the extraction of the scattering  operator 
from the autocorrelation of the measured data,
compared to the standard daylight imaging configuration addressed in
 \cite{claerbout68,claerbout85,noisebook,claerbout99},
in which the bottom condition is a radiating condition.
However, with arbitrary noise source distributions,
the autocorrelation function (\ref{eq:cor2a}) and the symmetrized response function (\ref{eq:symrep2}) have the same denominator, 
so that the identification of the eigenfrequencies can be carried out with both sets of data with 
the same accuracy.
Moreover, in a realistic source configuration the relationship 
  between the autocorrelation function (\ref{eq:cor2a}) and the scattering operator (\ref{eq:refop})
becomes much simpler.
\redc{In fact,  as we show in Eq. (\ref{eq:expl}),  the autocorrelation may produce  the 
scattering operator up to a frequency-dependent modulation function that depends
on the temporal spectrum of the noise source trace.  This is the case     
if we assume that the support of the noise sources (that is to say, the support of the function $K$)
is localized below the surface, and its thickness is smaller than the typical wavelength.
Then for any $R_s$ in the support of $K$, we have $\hat{\cal S}^\eps_{l}(\omega,R_s) \simeq 1$ by (\ref{eq:hatSeps}),
$\hat{\cal T}^\eps_{l}(\omega,R_s) \simeq \hat{\cal R}^\eps_{l}(\omega,R_o)$ by (\ref{def:Reps}) and (\ref{eq:hatTeps}), 
and $\tau(R_s,R_o) \simeq 0$, so the square brackets in (\ref{eq:cor2a}) can be simplified as
\begin{eqnarray*}
  \Big[ \int_0^{R_o}   \frac{c_o(R_s)}{R_s^2} K(R_s)  
\big|
\hat{\cal S}_{l}^\eps(\omega,R_s) 
+ \hat{\cal T}_{l}^\eps(\omega,R_s)
e^{-2 i \omega \tau(R_s,R_o)/\eps} 
\big|^2 dR_s \Big] \\
 \simeq 
\big|
1+ \hat{\cal R}_{l}^\eps(\omega,R_o)
\big|^2
 \Big[ \int_0^{R_o}   \frac{c_o(R_s)}{R_s^2} K(R_s)  
 dR_s \Big]  .
\end{eqnarray*}
As a consequence we get the following result.}
\begin{proposition}
\label{prop4}
If the spatial support of the noise source is localized in a small annulus below the surface, then
the statistical autocorrelation is of the form
\begin{eqnarray}
\label{eq:expl}
\fl
\hat{C}_{l,m}^\eps (\omega,R_o) 
\simeq \eps^2
 \hat{F}(\omega) \frac{c_o(R_o)}{2 R_o^2 \omega^2}   
\Big[ \int_0^{R_o}  \frac{c_o(R_s)}{R_s^2} K(R_s) dR_s
 \Big]
|1+\Gamma_{R_o}|^2 
\frac{\big| 1 + \hat{\cal R}_{l}^\eps(\omega,R_o) \big|^2 }
{\big| 1- \Gamma_{R_o} \hat{\cal R}_{l}^\eps(\omega,R_o)\big|^2} .
 \end{eqnarray}
 \end{proposition}
This result shows that the autocorrelation of the noise signals (\ref{eq:cor2a}) 
is directly related to the scattering operator (\ref{eq:refop}).
More precisely, the time derivative of the autocorrelation function is proportional to
the symmetrized response function (\ref{eq:responsesym}),
provided $\hat{F}(\omega)=\hat{f}(\omega)$ (which imposes in particular that the pulse profile $f$ should be an even function):
\begin{eqnarray}
\partial_t C_{l,m}(t)  \propto  p^{{\rm sym}}_{l,m}(t) .
\end{eqnarray}
We recover the classical seismic interferometry formula that has been established in many other configurations \cite{noisebook,schuster09,wap10a}.

 \redc{\section{On Robustness}\label{sec:rob} }

\subsection{Robustness with Respect to Sampling Configuration}
\label{rob1}
In the previous Section \ref{sec:correlation} 
we  studied the autocorrelation function $C_{l,m}(t,R_o)$.
Note that this is the autocorrelation function of the $(l,m)$-th component  associated with the spherical harmonics,
so that it is a linear transform of the cross correlation of the field observed at the surface:
\begin{eqnarray*}
\fl
C_{l,m}(t,R_o) &=& \int_0^\pi \int_0^{ 2\pi}  \int_0^\pi \int_0^{ 2\pi} \EE \big[ p(0,R_o,\theta,\varphi) 
 p(t,R_o,\theta',\varphi')\big]  Y_{l,m}(\theta,\varphi)  \overline{Y_{l,m}(\theta',\varphi')} \\
\fl && \times 
\sin \theta d\varphi d \theta
\sin \theta'd\varphi' d \theta' .
\end{eqnarray*}
From this it may seem that we need to 
observe the field everywhere at the surface to estimate this quantity (by further 
substituting  a  time average for the expectation).
However, we do not need such extensive data \redc{and we can relax the hypothesis that we 
observe the field everywhere at the surface of the Earth.}
Indeed we have shown that $C_{l,m}(t,R_o) $ does not depend on $m$.
This is due to the spherical symmetry of the Earth model and the  noise source distribution.
Therefore we can consider
\begin{eqnarray*}
C_l(t,R_o) = \frac{1}{2l+1} \sum_{m=-l}^l C_{l,m}(t,R_o) 
\end{eqnarray*}
where actually all terms are equal and given by the expression in Proposition \ref{prop:statauto}.
Using the addition theorem of spherical harmonics,
we have
\begin{eqnarray*}
C_l(t,R_o)&=&\int_0^\pi \int_0^{ 2\pi}  \int_0^\pi \int_0^{ 2\pi} \EE \big[ p(0,R_o,\theta,\varphi) 
 p(t,R_o,\theta',\varphi')\big] \frac{1}{4\pi} P_l(\cos \Omega) \\
&& \times \sin \theta d\varphi d \theta
\sin \theta'd\varphi' d \theta' ,
\end{eqnarray*}
where $P_l$ is the Legendre polynomial and $\cos \Omega$ is the angle between two unit vectors oriented 
 at the polar coordinates $(\theta,\varphi)$ and $(\theta',\varphi')$.
The expectation 
\begin{eqnarray*}
{\cal C}(t,\Omega) = \EE \big[ p(0,R_o,\theta,\varphi) 
 p(t,R_o,\theta',\varphi')\big] 
\end{eqnarray*}
 depends only on $\cos \Omega$ and it can be expanded as
\begin{eqnarray*}
 {\cal C}(t,\Omega) = \sum_{l=0}^\infty C_l(t,R_o) P_l(\cos \Omega) .
\end{eqnarray*}
The quantity of interest $C_l(t,R_o)$ is given by
\begin{eqnarray}
\label{eq:expressCl}
C_l(t,R_o) = \frac{2l+1}{2}  \int_0^\pi {\cal C}(t,\Omega) P_l(\cos \Omega) \sin \Omega d\Omega .
\end{eqnarray}
 \redc{This shows that observation points at the surface of the Earth that are such that the angles $\Omega$ between
 pairs of points cover the interval $(0,\pi)$ are sufficient to estimate the Earth spectrum (and the scattering operator
 under the additional hypothesis in Proposition \ref{prop4}): 
 this data set gives estimates of $ {\cal C}(t,\Omega)$
 for a sufficiently dense grid of $\Omega$; then the quantity of interest $C_l(t,R_o)$ can be obtained 
 by a numerical evaluation (a quadrature formula) of the integral (\ref{eq:expressCl});
the Earth point spectrum can then be obtained by inspection of the Fourier transform $\hat{C}_l(\omega,R_o)$
whose peaks correspond to eigenfrequencies.
%  This important result indeed shows how data recordings should be effectively combined so that 
% they can be related to the Earth point spectrum.
 }

\subsection{Robustness with Respect to Source Configuration}
\label{rob2} %
Consider the case when the assumption in (\ref{eq:conf})
is generalized as 
\begin{eqnarray*}
\fl
\EE \big[ {f}(t,r,\theta,\varphi){f}(t',r',\theta',\varphi') ] = F (  {t-t'}  )
K(r)G(\theta,\varphi) \delta(r-r') \sin(\theta)^{-1} \delta(\theta-\theta') \delta(\varphi-\varphi') ,
%\label{eq:conf2}
\end{eqnarray*}
thus allowing for general lateral source distribution with density function $G(\theta,\varphi)$. 
Then 
\begin{eqnarray*}
%\label{eq:spec2}
\EE \big[ \hat{f}_{l,m}(\omega,r) \overline{\hat{f}_{l',m'}(\omega',r')} \big] = 2\pi
\hat{F}(\omega) K(r)  G_{l,m} \delta(r-r') \delta(\omega-\omega') \delta_{ll'}\delta_{mm'} ,
\end{eqnarray*}
with
$$
G_{l,m} =   \int_0^\pi \int_0^{ 2\pi} |Y_{lm}(\theta,\varphi) |^2 G(\theta,\varphi) \sin \theta  d\varphi d\theta.
$$
In this general setting we find that Propositions \ref{prop:statauto} and \ref{prop4} still hold true but $\hat{C}_{l,m}^\eps (\omega,R_o)$ is
now multiplied by $G_{l,m}$ compared to Eqs.~(\ref{eq:cor2a}) and (\ref{eq:expl}).
%\begin{eqnarray*}
%\hat{C}_{l,m}^\eps (\omega,R_o) \mapsto  G_{l,m}  \hat{C}_{l,m}^\eps (\omega,R_o)   ,
%\end{eqnarray*}
%with  $\hat{C}_{l,m}^\eps$ defined as before and in particular being independent of $m$.
Since these new expressions depend on $m$ only though $G_{l,m}$, 
% right-hand sides of Eqs.~(\ref{eq:cor2a}) and (\ref{eq:expl}) 
we can sum over $m=-l,\ldots,l$ to get \redc{rid of} the dependence with respect to $m$
of the source distribution, and \redc{there remains} only a dependence with respect to $l$ via a positive multiplicative
factor $\sum_{m=-l}^l G_{l,m}$.
This shows that, as long as the main objective is to extract the eigenfrequencies,
we have robustness with respect 
to the assumption of angular homogeneity in the random source distribution.
However, the quantitative estimation of the Earth scattering operator is sensitive 
to the angular source distribution through the multiplicative factor  $\sum_{m=-l}^l G_{l,m}$
that affects the amplitude of the estimation.

\subsection{Robutsness with Respect to Medium Noise or Model Error}
\label{rob3}%
We have assumed in the previous sections that the speed of propagation has only radial variations.
Here we want to show that the result is robust with respect to certain angular variations of the velocity model. 
More exactly, we consider an Earth-model with slow angular variations, whose velocity model has the form:
\begin{eqnarray}
\frac{1}{c^2(r,\theta,\varphi)} = \frac{1}{c_o^2(r)} + V_1^\eps(r)
+
V_2^\eps(r,\theta,\varphi) ,
\end{eqnarray}
in which there are small and slow angular velocity fluctuations $V_2^\eps$. 
We call them slow because they vary at the scale one with respect to the angles $\theta$ and $\varphi$, 
while they may vary slowly or rapidly in the radial coordinate $r$.
We want to clarify under which circumstances the perturbation  $V_2^\eps$ can be neglected.

The field has the form (\ref{eq:fieldmode1a}-\ref{eq:fieldmode1b})
and the coupled system of first-order ODEs for the mode amplitudes now reads:
\begin{eqnarray}
\nonumber
\fl
\frac{\partial \hat{\alpha}_{l,m}^\eps}{\partial r} (\omega,r)&=&
- i \frac{\omega}{2 \eps} c_o(r)V_1^\eps(r) 
 \big[   \hat{\alpha}_{l,m}^\eps(\omega,r)+ e^{-2i \frac{\omega\tau(r,R_o)}{\eps}} \hat{\beta}_{l,m}^\eps(\omega,r)\big]\\
\label{sys:4a}
\fl &&- i \frac{\omega}{2} c_o(r) 
\sum_{l',m'} q_{l,m,l',m'}^\eps(r)
 \big[   \hat{\alpha}_{l',m'}^\eps(\omega,r)+ e^{-2i \frac{\omega\tau(r,R_o)}{\eps}} \hat{\beta}_{l',m'}^\eps(\omega,r)\big] ,\\
 \nonumber
\fl \frac{\partial \hat{\beta}_{l,m}^\eps}{\partial r} (\omega,r)&=&
 i \frac{\omega}{2 \eps} c_o(r)  V_1^\eps(r) 
  \big[ e^{2i \frac{\omega\tau(r,R_o)}{\eps}}  \hat{\alpha}_{l,m}^\eps(\omega,r)+   \hat{\beta}_{l,m}^\eps(\omega,r)\big] \\
\fl && + i \frac{\omega }{2 } c_o(r)  
   \sum_{l',m'} q_{l,m,l',m'}^\eps(r)
 \big[ e^{2i \frac{\omega\tau(r,R_o)}{\eps}}  \hat{\alpha}_{l',m'}^\eps(\omega,r)+   \hat{\beta}_{l',m'}^\eps(\omega,r)\big] ,
 \label{sys:4b}
\end{eqnarray}
for $r \in (0,R_o)$, instead of (\ref{sys:3a}-\ref{sys:3b}),
with
 \begin{eqnarray*}
q_{l,m,l',m'}^\eps (r) &=&  \frac{1}{\eps}
\int_0^\pi \int_0^{ 2\pi}  \overline{Y_{l,m}(\theta,\varphi) }
V^\eps_2(r,\theta,\varphi) Y_{l',m'}(\theta,\varphi) \sin \theta d\varphi d \theta .
\end{eqnarray*}
We consider the case in which the fluctuation term $V_2^\eps$ may have slow and fast components:
\begin{eqnarray}
V_2^\eps (r,\theta,\varphi) = \eps^a V_{21}(r,\theta,\varphi) 
+\eps^b V_{22} \big( \frac{r}{\eps},\theta,\varphi \big) 
+\eps^c V_{23} \big( \frac{r}{\eps^2},\theta,\varphi \big) ,
\end{eqnarray}
where 
$V_{21}$ and $V_{22}$ can be \redc{arbitrary functions} and 
$V_{23}$ is a zero-mean random process, that is stationary and mixing in $r$.

We have the following results:\\
- Provided $a>1$, the term $V_{21}$ gives rise to terms in  (\ref{sys:4a}-\ref{sys:4b}) 
that vanish in the limit $\eps \to 0$.\\
-
Provided $b>1$, the term $V_{22}$  gives rise to terms in  (\ref{sys:4a}-\ref{sys:4b}) 
that vanish in the limit $\eps \to 0$.\\
- 
Provided $c>0$, the term $V_{23}$  gives rise to terms in  (\ref{sys:4a}-\ref{sys:4b}) 
that vanish in the limit $\eps \to 0$.

The first two assertions are trivial by direct inspection of the amplitude of the coupling terms $\eps^{a-1} q(r)$
and $\eps^{b-1} q(r/\eps)$, respectively,
but the third assertion is not so trivial as it involves coupling terms of the form $\eps^{c-1} q(r/\eps^2)$
\redc{for}  a mixing and zero-mean process $q$.
However, diffusion approximation theory reveals that this coupling becomes effective only when $c=0$ \cite{book}.
Therefore all crossed terms in the equations (\ref{sys:4a}-\ref{sys:4b}) cancel  
and we get a system similar to (\ref{sys:3a}-\ref{sys:3b}).
The energy conservation equation (\ref{eq:ener}), 
which is the key property, holds true and
the conclusion of the previous section holds  including robustness with respect to
partial source and measurement aperture.

\subsection{Robustness with Respect to Measurement Noise}
\label{rob4}%
The results presented in the previous sections are based on the behavior of the statistical autocorrelation function $C_{l,m}(t,R_o)$.
In this section we want to explain that this function can be estimated by the empirical autocorrelation of the recorded
signals and that it is quite robust with respect to additive measurement noise.

The statistical autocorrelation function  $C_{l,m}(t,R_o)$ defined by (\ref{def:corstat}) is the one that is related to the scattering operator.
The empirical autocorrelation function  $C_{l,m}^T(t,R_o)$ defined by (\ref{eq:empiricalcov}) is the one that is computed
from the recorded data. 
It is easy to check that the expectation of $C_{l,m}^T(t,R_o)$ is exactly $C_{l,m}(t,R_o)$ for any $T$.
By a detailed fluctuation analysis it is also possible to show that the variance of $C_{l,m}^T(t,R_o)$
is proportional to $1/T$ \cite[Chapter 2]{noisebook}. This ensures that the empirical autocorrelation
function is a good estimate of the statistical autocorrelation function provided the recording time window 
$T$ is large enough. As noticed in \cite[Chapter 2]{noisebook}, the required recording time is all the longer as 
lower frequency components are being investigated.

When the recorded signals are polluted by additive measurement noise,
which is independent from the noise sources,
the empirical autocorrelation function in Eq. 
(\ref{eq:empiricalcov}),   computed via a discretization of
Eq.  (\ref{eq:empiricalcov2}),   is the sum of the statistical autocorrelation function  $C_{l,m}(t,R_o)$ 
and of the autocorrelation function $C_{l,m}^{\rm n}(t)$ that comes from the additive noise
at the $N$ receivers.
Measurement noise is independent from one receiver to the other one, so that, if the $N$ receivers
are distributed uniformly at the surface of the Earth and the measurement noise statistics is the same at each receiver we have:
$$
C_{l,m}^{\rm n}(t) =\frac{1}{N} F^{\rm n}(t) , 
$$
where $F^{\rm n}(t)$ is the autocorrelation function of the measurement noise ($\hat{F}^{\rm n}(\omega)$
is the power spectral density of the measurement noise). 
This shows that the impact of the measurement noise decays with the number of receivers.
Moreover, measurement noise usually has higher frequencies than the spectral band that is of interest
for global seismology.
However, if the spectrum of the additive noise intersects the spectral band over which we look for eigenfrequencies
and/or the scattering operator, then measurement noise may corrupt the estimation.
The estimation of the eigenfrequencies should be quite robust as measurement noise is very unlikely to
have spectral peaks similar to the ones that are investigated.
The estimation of the scattering operator (which is more sensitive to the estimated amplitude)
may be affected.

\section{Concluding Remarks}

We have shown  that global acoustic daylight imaging is 
possible. 
The point spectrum of the Earth can be extracted from  the correlation functions of the signals recorded at the surface and emitted by 
unknown noise sources localized away from the surface.
Under an additional realistic assumption of the spatial support of the noise sources,
the complete scattering operator can be extracted from the correlation functions.
\redc{The first result (on the extraction of the point spectrum)
is robust to receiver distribution, source distribution, and medium and measurement noise, 
but one should be careful to interpret the second result when the noise sources are not evenly distributed 
as this may perturb the amplitudes of the recovered signals.}

In this paper we have only considered low angular modes (with Legendre number $l$ of the order of one). 
Indeed we do not need to analyze the high angular modes because there is no coupling between low and high angular modes in our setting. 
We could consider high angular modes, but then a non-trapping condition would be necessary, such as the Herglotz \cite{herglotz} and Wiechert and Zoeppritz \cite{wiechert} conditions.

We have here,  in Eq. (\ref{eq:source2}), considered the situation 
 that the noise sources are delta-correlated in space. 
This assumption simplifies the analysis but can be relaxed provided the correlation radius is small 
by using stationary phase methods~\cite{noisebook} or 
semi-classical analysis~\cite{colin}.

\section{Acknowledgement}
Maarten V. de Hoop was supported in part by the Simons Foundation under the MATH + X program.  Knut  S\o lna was supported in part by  
AFOSR under grant FA9550-14-1-0197 and NSF under grant 1616954.

\end{document}